\documentclass[10pt,leqno]{amsart}
\usepackage{graphicx}
\usepackage{indentfirst,csquotes}
\usepackage{hyperref}
\usepackage{amssymb,amsthm,amsmath, titlesec}
\usepackage{xcolor,paralist,hyperref,titlesec,fancyhdr,etoolbox}
\newtheorem{theorem}{Theorem}[]

\newtheorem{example}[theorem]{Example}
\newtheorem{lemma}[theorem]{Lemma}
\newtheorem{proposition}[theorem]{Proposition}
\newtheorem{corollary}[theorem]{Corollary}

\newtheorem{remark}[theorem]{Remark}

\titleformat{\section}[display]{\normalfont\huge\centering}{\centering\chaptertitlename\thechapter}{10pt}{\large}
\titlespacing*{\section}{0pt}{0ex}{0ex}

\hypersetup{ colorlinks=true, linkcolor=black, filecolor=black, urlcolor=black }

\usepackage{lipsum}

\begin{document}
\title{Prime and weakly prime submodules on amalgamated duplication of a ring along an ideal} 
\author[GURSEL YESILOT, ESRA TARAKCI AND YASEMIN SIMSEK]{GURSEL YESILOT, ESRA TARAKCI AND YASEMIN SIMSEK}
\date{\today}
\address{Yildiz Technical University, Department of Mathematics, Davutpaşa, Istanbul, Turkey}
\email{gyesilot@yildiz.edu.tr}
\address{Dogus University, Department of Computer Engineering, Istanbul, Turkey}
\email{etarakci@dogus.edu.tr}
\address{Maltepe University, Department of Mathematics and Science Education, Istanbul, Turkey}
\email{yaseminsimsek@maltepe.edu.tr}
\maketitle

\let\thefootnote\relax

\begin{abstract}
Let $A$ be a commutative ring with identity. A proper submodule $N$ of $A$-module $M$ is said to be prime submodule if $ax \in N$ where $a \in A, x \in M$, implies $x \in N$ or $aM \subseteq N$. A proper submodule $N \subset M$ is said to be weakly prime submodule if $0 \neq ax \in N$ where $a \in A, x \in M$, then either $x \in N$ or $aM \subseteq N$. The notion of weakly prime submodule was introduced by Atani and Farzalipour \cite{atani2007weakly}. The purpose of this paper is to study the form of prime and weakly prime submodules of duplication of the $A$-module $M$ along the ideal $I$ (denoted by $M \bowtie I$), introduced and studied by E. M. Bouba, N. Mahdou and M. Tamekkante. A number of results concerning prime and weakly prime submodules on amalgamated duplication and examples are given. 
\end{abstract} 

\bigskip
\begin{center}
\title{1. INTRODUCTION}    
\end{center}

In this article, we assume all rings are commutative with identity and all modules are unital. We will start with certain definitions and notations that will be used throughout the paper.
 A proper submodule $N$ of $A$-module $M$ is called prime submodule if $ax \in N$ where $a \in A$ and $x \in M$, then either $x \in N$ or $aM \subseteq N$.
 A proper submodule $N$ of $A$-module $M$ is called weakly prime submodule if $0 \neq ax \in N$ where $a \in A, x \in M$ implies
$x \in N$ or $aM \subseteq N$. In 2007, S. E. Atani ve F. Farzalipour  aimed to prove some results for weakly prime submodules in their article "On Weakly Prime Submodules" \cite{atani2007weakly}.
Equivalent to the above definition, A. Azizi defined the weakly prime submodule as "Let $A$ be a commutative ring with identity and $M$ be a $A$-module. For $a,b \in A$ and every submodule $T$ of $M$, if $abT \subseteq N$ implies $aT \subseteq N$ or $bT \subseteq N$, then $N$ is called a weakly prime submodule". \\
Recall that a proper submodule $N$ of $M$ is said to be a primary submodule if the condition $ax \in N$, $a \in A$ and $x \in M$, implies that $x \in N$ or $a^nM \subseteq N$, for some
positive number $n$. 
 
 In 2007, the construction $A \bowtie I$ was introduced and its basic properties were studied by D’Anna and Fontana in \cite{d2007amalgamated}.\\
 Let $A$ be a commutative ring and $I$ be an ideal of $A$. The amalgamated duplication of $A$ along $I$, denoted by $A \bowtie I$, is the special subring of $A \times A$ defined by $A \bowtie I := \{(a, a + i) | a \in A, i \in I\}$.\\
 $A$ is embedded in $A \bowtie I$ by the mapping $a \mapsto (a,a)$, and this new ring can be thought of as an extension of $A$.\\
 In 2018, E. M. Bouba, N. Mahdou and M. Tamekkante studied on the structure of the amalgamated algebra with modules in their article.
 Let $I$ be an ideal of $A$, and $M$ an $A$-module. The duplication of the $A$-module $M$ along the ideal $I$ denoted by $M \bowtie I$ and gave the following definition 
\begin{center}
    $M \bowtie I := \{(m,m') \in M \times M| m - m' \in IM\}$
\end{center}
 which is an $A \bowtie I$-module with the multiplication given by $(a, a+ i).(m, m') = (am, (a+i)m')$, where $a \in A, i \in I$, and $(m,m') \in M \bowtie I$. \\
 It is clear that if $M = A$, then the duplication of the $A$-module $A$ along the ideal $I$ coincides with the amalgamated duplication of the ring $A$ along the ideal $I$.\\
 
 It is defined as $N \bowtie I := \{(n,m) \in N \times M, n-m \in IM\}$, where $N$ is a proper submodule of $M$. Therefore $N \bowtie I$ is a submodule of $M \bowtie I$.\\
In this paper, we begin a study of the basic properties of prime and weakly prime submodules of $M \bowtie I$.\\
In this paper, we pursue the investigation on the structure of the form
$M \bowtie I$, with a particular attention to the form of prime submodules and to weakly prime submodules introduced and studied \cite{atani2007weakly}.
\\
\,\\
\begin{center}
\title\MakeUppercase{2. Properties of prime submodules, weakly prime submodules and primary submodules on amalgamated duplication}    
\end{center}

To avoid unnecessary repetition, let us fix notation for the rest of the paper.
All along this paper, $A \bowtie I$ will denote the amalgamated duplication of $A$ along $I$ and $M \bowtie I$ will denote the duplication of the $A$-module $M$ along the ideal $I$.
$0 \times I$ is an ideal of $A \bowtie I$ and also,
$(0 \times I)M \bowtie I = \{(0, i)(m_1, m_2)|i \in I,(m_1, m_2) \in M \times M, m_1 - m_2 \in IM \} = 0 \times IM$ and $IM \times IM$ are $A \bowtie I$-submodules of $M \bowtie I$. Moreover, let $N$ be the submodule of $M$.
It is defined as $N \bowtie I = \{(n_1, m_1) \in N \times M|n_1 - m_1 \in IM\}$. Therefore $N \bowtie I$ is a submodule of $M \bowtie I$. Our first result is a characterization of prime submodules of the form $N \bowtie I$ of duplication of $A$-module $M$ along the ideal $I$.

\begin{lemma}
\label{Lemma 1.}
     Let $N$ be a submodule of $M$. Then $(N \bowtie I : M \bowtie I) = (N : M) \bowtie I$ \cite{issoual2022s}.
\end{lemma}

\begin{lemma}
    Let $L$ be a submodule of $M$. Then the following statements are equivalent \cite{issoual2022s}:
    \begin{enumerate}
        \item  $L$ is a prime submodule of $M$.
        \item  $L \bowtie I$ is a prime submodule of $M \bowtie I$.
    \end{enumerate}
\end{lemma}

\begin{proof}
    Follows from \cite{issoual2022s}.
\end{proof}

\begin{example}
   Let $A = \mathbb{Z}$,  $\mathbb{Z}$-module $M = \mathbb{Z}$ and $L = 3\mathbb{Z}$. Since every prime ideal of $A$ is a prime submodule of the $A$-module $A$, $L$ is a prime submodule of $M$. For $I = 4\mathbb{Z}$,  $L \bowtie I = 3\mathbb{Z} \bowtie 4\mathbb{Z}$ is a prime submodule of $M \bowtie I =\mathbb{Z} \bowtie 4\mathbb{Z} $. If we generalize, for $I = n\mathbb{Z}, n > 1$, we conclude that $L \bowtie I$ is a prime submodule.
\end{example}

\begin{corollary}
     Let $L$ be a submodule of $M$. Then the following statements are equivalent:
    \begin{enumerate}
        \item  $L$ is a weakly prime submodule of $M$.
        \item  $L \bowtie I$ is a weakly prime submodule of $M \bowtie I$.
    \end{enumerate}
\end{corollary}

\begin{proof}
   $(1) \Rightarrow (2)$ Assume that $L$ is a weakly prime submodule of $M$. \\Let $(0,0) \neq (a,a+i)(m,m') \in L \bowtie I$ where $(a,a+i) \in A \bowtie I$ and $(m,m') \in M \bowtie I$. Hence, $(0,0) \neq (am, (a+i)m') \in L \bowtie I$ and we observe that $0 \neq am \in L$. Since $L$ is a weakly prime submodule, we conclude either $m \in L$ or $a \in (L:M)$. \\
   If $m \in L $, since $m-m' \in IM$, then $(m,m') \in L \bowtie I$.\\
   If $a  \in (L: M)$, according to Lemma \ref{Lemma 1.}, then we have $(a, a +i)  \in (L:M) \bowtie I = (L \bowtie I : M \bowtie I) $. Therefore, $L \bowtie I$ is a weakly prime submodule.\\
   $(2) \Rightarrow (1)$ Suppose that $L \bowtie I$ is a weakly prime submodule. Let $0 \neq a_1m \in L$ where $a_1 \in A$ and $m \in M$. Hence, we have $(0,0) \neq (a_1, a_1)(m,m) \in L \bowtie I$. Since $L \bowtie I $ is a weakly prime submodule, we have either $(m,m) \in L \bowtie I$ or $(a_1,a_1) \in (L \bowtie I : M \bowtie I)$. \\
   If $(m,m) \in L \bowtie I$, then $m \in L$.\\
   If  $(a_1,a_1) \in (L \bowtie I : M \bowtie I)$, according to Lemma \ref{Lemma 1.}, we have $(a_1,a_1) \in (L: M) \bowtie I$ and so $a_1 \in (L:M)$. Therefore, $L$ is a weakly prime submodule of $M$.
   
\end{proof}

\begin{example}
    Let $A = \mathbb{Z}_6$, $\mathbb{Z}_6$-module $M = \mathbb{Z}_6$ and $L = \{\Bar{0}\}$. $\{\Bar{0}\}$ is a weakly prime submodule. For $I = 3\mathbb{Z}_6$, $L \bowtie I = \{\Bar{0}\} \bowtie 3\mathbb{Z}_6$ is a weakly prime submodule of $M \bowtie I = \mathbb{Z}_6 \bowtie 3\mathbb{Z}_6$.
\end{example}

As it can be understood from the definitions of prime and weakly prime submodule, every prime submodule is a weakly prime. In contrast, weakly prime submodule need not be prime as in the following example.

\begin{example}
     For $\mathbb{Z}_n$-module $\mathbb{Z}_n$, where $n$ is a composite number, $0$ is a weakly prime submodule but not a prime. For instance, let $A = \mathbb{Z}_6$, $\mathbb{Z}_6$-module $M = \mathbb{Z}_6$, $I =3\mathbb{Z}_6$ and $L = \{\Bar{0}\}$. According to the previous example, $L \bowtie I$ is a weakly prime submodule, however, while $(\Bar{3}, \Bar{3})(\Bar{4}, \Bar{1}) = (\Bar{0}, \Bar{3}) \in L \bowtie I$, we have $(\Bar{4}, \Bar{1}) \notin L \bowtie I$ and $(\Bar{3}, \Bar{3}) M \bowtie I \nsubseteq L \bowtie I$.
\end{example}

\begin{example}
 $L$ is a primary submodule of $M$ if and only if $L \bowtie I$ is a primary submodule of $M \bowtie I$. Let consider the following example.
  Let $A = \mathbb{Z}$, $\mathbb{Z}$-module $M = \mathbb{Z}$ and $L = 5\mathbb{Z}$ be submodule of $M$. Since every prime ideal is a primary and also every primary ideal in $A$ is a primary submodule of the $A$-module $A$, $L$ is a primary submodule. For $I = 4\mathbb{Z}$, $L \bowtie I = 5\mathbb{Z} \bowtie 4\mathbb{Z} $ is a primary submodule of $M \bowtie I = \mathbb{Z} \bowtie 4\mathbb{Z}$. As the following proposition concludes that this example is also true in general.
\end{example}

\begin{proposition}
    Let $L$ be a submodule of $M$. Then the following statements are equivalent:
    \begin{enumerate}
        \item  $L$ is a primary submodule of $M$.
        \item  $L \bowtie I$ is a primary submodule of $M \bowtie I$.
    \end{enumerate}
\end{proposition}

\begin{proof}
    $(1) \Rightarrow (2)$ Assume that $L$ is a primary submodule of $M$.  Let $(a,a+i)(x,x') \in L \bowtie I$ and $(x,x') \notin L \bowtie I$ where $(a,a+i) \in A \bowtie I$, $(x,x') \in M \bowtie I$. Hence, we have $(ax, (a+i)x')) \in L \bowtie I$ and $ax \in L$. Since $x \notin L$ and $L$ is a primary submodule, then we have $a^n \in (L :M)$. Therefore, $(a, a+i)^n = (a^n, a^n + i') \in (L:M) \bowtie I = (L \bowtie I: M \bowtie I)$ where $i' = \sum_{k=0}^{n-1} \binom{n}{k} a^k i^{n-k} $. So, we conclude that $ (a, a+i)^n \in (L \bowtie I: M \bowtie I)$ and $L \bowtie I$ is a primary submodule of $M \bowtie I$.\\
    $(2) \Rightarrow (1)$ Suppose that $L \bowtie I$ is a primary submodule. Now let $ax \in L$ and for all $t \in \mathbb{Z}^{+}$,  $a^t \notin (L:M)$ where $a \in A, x \in M$. Hence we have $(a,a)(x,x) \in L \bowtie I$. Since $a^t \notin (L: M)$, then $(a,a)^t = (a^t, a^t)  \notin (L:M) \bowtie I = (L \bowtie I : M \bowtie I)$. Since $(a,a)^t \notin (L \bowtie I : M \bowtie I) $ and $L \bowtie I$ is a primary, then we conclude $(x,x) \in L \bowtie I$. So, $x \in L$ and L is a primary submodule of $A$.
\end{proof}

As it can be understood from the definition, every prime submodule is a primary. However, primary submodule need not be a prime as in the following example.
\begin{example}
     Let $A = \mathbb{Z}$, $\mathbb{Z}$-module $M = \mathbb{Z}$ and $L = 8\mathbb{Z}$ be submodule of $M$. For $I = 4\mathbb{Z}$, $L \bowtie I = 8\mathbb{Z} \bowtie 4\mathbb{Z} $ is a primary submodule of $\mathbb{Z} \bowtie  4\mathbb{Z}$. However, while $(2,6)(4,8) = (8,48) \in L \bowtie I$, we observe $(4,8) \notin L \bowtie I$ and $(2,6)M \bowtie I \nsubseteq L \bowtie I$. Therefore  $L \bowtie I$ is not a prime submodule. 
\end{example}

\begin{lemma} 
\label{Lemma 3.}
  Let $M \bowtie I$ be an $A \bowtie I$-module and $N \bowtie I$ a proper submodule of $M \bowtie I$.
  \begin{itemize}
      \item[(i)]  $N \bowtie I$ is a weakly prime submodule if and only if for every submodule $K \bowtie I$ of $M \bowtie I$ not contained in $N \bowtie I$, $(N \bowtie I : K \bowtie I)$ is a prime ideal of $A \bowtie I$. In particular $(N \bowtie I: M \bowtie I)$ is a prime ideal of $A \bowtie I$.
      \item[(ii)]  Let $N \bowtie I$ be a weakly prime submodule of $M \bowtie I$. Then for all submodules $K \bowtie I$ and $L \bowtie I$ of $M \bowtie I$ not contained in $N \bowtie I$, $(N \bowtie I: K \bowtie I) \subseteq (N \bowtie I : L \bowtie I)$ or $(N \bowtie I : L \bowtie I) \subseteq (N \bowtie I : K \bowtie I)$.
  \end{itemize}

\end{lemma}

\begin{proof}
    \begin{itemize}
        \item[(i)] $(\Rightarrow):$  Let $N \bowtie I$ be a weakly prime submodule and $K \bowtie I \nsubseteq N \bowtie I$. Then, $(N \bowtie I: K \bowtie I) \neq A \bowtie I$. If $(N \bowtie I : K \bowtie I) = A \bowtie I$, there would be $(a,a+i)K \bowtie I \subseteq N \bowtie I$  for every $(a,a+i) \in A \bowtie I$. Also it would be a contradiction that $K \bowtie I = (1,1)K \bowtie I \subseteq N \bowtie I$ for $(1,1) \in A \bowtie I$. So, $(N \bowtie I : K \bowtie I) $ is a proper ideal.\\
        Now let take $(a_1,a_1+i_1), (a_2,a_2+i_2) \in A \bowtie I$ such that $(a_1,a_1+i_1)(a_2,a_2+i_2) \in (N \bowtie I: K \bowtie I)$. Then we have $(a_1,a_1+i_1)(a_2,a_2+i_2)K \bowtie I \subseteq N \bowtie I$. Since $N \bowtie I$ is a weakly prime submodule, then $(a_1,a_1+i_1)K \bowtie I \subseteq N \bowtie I$ or $(a_2,a_2+i_2)K \bowtie I \subseteq N \bowtie I$. \\
        If $(a_1,a_1+i_1)K \bowtie I \subseteq N \bowtie I$, then $(a_1,a_1+i_1) \in (N \bowtie I: K \bowtie I)$.\\
        If $(a_2,a_2+i_2)K \bowtie I \subseteq N \bowtie I$,then  $(a_2,a_2+i_2) \in (N \bowtie I:K \bowtie I)$.\\
        So, $(N \bowtie I : K \bowtie I)$ is a prime ideal of $A \bowtie I$. \\
        
        In particular, since $N \bowtie I$ is a proper submodule, $M \bowtie I$ is a submodule of itself that is not contained in $N \bowtie I$, and it follows that $(N \bowtie I: M \bowtie I)$ is also a prime ideal of $A \bowtie I$.\\

        $(\Leftarrow):$ Let $(N \bowtie I:K \bowtie I)$ be a prime ideal of $A \bowtie I$. Then $(N \bowtie I: K \bowtie I) \neq A \bowtie I$, that is, there is an $(a,a+i) \in A \bowtie I$, for which $(a,a+i)K \bowtie I \nsubseteq N \bowtie I$. \\
        If $K \bowtie I \subseteq N \bowtie I$, then $(a,a+i)K \bowtie I \subseteq K \bowtie I \subseteq N \bowtie I$ and there would be a contradiction $(a,a+i)K \bowtie I \subseteq N \bowtie I$. Therefore $K \bowtie I \nsubseteq N \bowtie I$. Now let $(a_1,a_1 + i_1)(a_2, a_2 + i_2)K \bowtie I \subseteq N \bowtie I$ for $(a_1,a_1 + i_1), (a_2, a_2 + i_2) \in A \bowtie I$. Hence $ (a_1,a_1 + i_1)(a_2, a_2 + i_2) \in (N \bowtie I: K \bowtie I)$. Since $(N \bowtie I: K \bowtie I)$ is a prime ideal, then either $(a_1,a_1 + i_1) \in (N \bowtie I: K \bowtie I)$ or $(a_2, a_2 + i_2) \in (N \bowtie I: K \bowtie I)$. That is, $(a_1,a_1 + i_1)K \bowtie I \subseteq N \bowtie I$ or $(a_2,a_2 + i_2) K \bowtie I \subseteq N \bowtie I$ is obtained, and accordingly $N \bowtie I$ is a weakly prime submodule. \\

        \item[(ii)] Let $N \bowtie I$ be a weakly prime submodule of $M \bowtie I$ and $K \bowtie I$, $L \bowtie I \nsubseteq N \bowtie I$. Assume that $(N \bowtie I: K \bowtie I) \nsubseteq (N \bowtie I : L \bowtie I)$ and $(N \bowtie I : L \bowtie I) \nsubseteq (N \bowtie I : K \bowtie I)$.\\
        That is, there is at least one element $(0,0) \neq (a_1,a_1+i_1) \in (N \bowtie I : K \bowtie I) \setminus (N \bowtie I: L \bowtie I)$ and $(0,0) \neq (a_2,a_2 +i _2) \in (N \bowtie I : L \bowtie I) \setminus (N \bowtie I : K \bowtie I)$.\\
        Hence we have $(a_1,a_1 + i_1)K \bowtie I \subseteq N \bowtie I$ and $(a_2,a_2 +i_2)L \bowtie I \subseteq N \bowtie I$. For every $(k,k') \in K \bowtie I$,  $(a_1, a_1+i_1)(k,k') \in  N \bowtie I$ can be written. \\
        Therefore, for every $(a,a+i) \in A \bowtie I, (a,a+i)(a_1,a_1+i_1)(k,k') \in N \bowtie I$, and specifically for $(a_2,a_2+i_2), \\
        (a_2,a_2+i_2)(a_1,a_1+i_1)(k,k') \in N \bowtie I $, $(\forall (k,k') \in K \bowtie I)$. Due to the fact that $K \bowtie I \nsubseteq N \bowtie I$ and $N \bowtie I$ is weakly prime submodule, there is at least one $(0,0) \neq (t,t') \in K \bowtie I \setminus N \bowtie I $ and for $(t,t')$, $(a_2,a_2+i_2)(a_1,a_1+i_1)M \bowtie I \subseteq N \bowtie I$ while $(a_2,a_2+i_2)(a_1,a_1+i_1)(t,t') \in N \bowtie I$.\\
        
        Therefore we conclude $(a_2,a_2+i_2)(a_1,a_1+i_1) \in (N \bowtie I : M \bowtie I)$. Since $(N \bowtie I : M \bowtie I)$ is prime ideal by $(i)$, we have $(a_2,a_2+i_2) \in (N \bowtie I : M \bowtie I)$ or $(a_1,a_1+i_1) \in (N \bowtie I : M \bowtie I)$.\\
        If $(a_2,a_2+i_2) \in (N \bowtie I: M \bowtie I)$, since $(N \bowtie I : M \bowtie I) \subseteq (N \bowtie I : K \bowtie I)$, the contradiction $(a_2,a_2+i_2) \in (N \bowtie I: K \bowtie I)$ is obtained.\\
        If $(a_1,a_1+i_1) \in (N \bowtie I : M \bowtie I)$, since $(N \bowtie I : M \bowtie I) \subseteq (N \bowtie I : L \bowtie I)$, then the contradiction $(a_1,a_1+i_1) \in (N \bowtie I: L \bowtie I)$ is obtained.\\
        Similarly, $(a_2,a_2+i_2)L \bowtie I \subseteq N \bowtie I$ leads us to the contradictions $(a_1,a_1+i_1) \in (N \bowtie I: L \bowtie I)$ or $(a_2, a_2+i_2) \in (N \bowtie I: K \bowtie I)$. Then it must be $(N \bowtie I : K \bowtie I) \subseteq (N \bowtie I: L \bowtie I)$ or $(N \bowtie I: L \bowtie I) \subseteq (N \bowtie I: K \bowtie I)$.

    \end{itemize}
\end{proof}

\begin{corollary} 
    Let $M \bowtie I$ be an $A \bowtie I$-module and $N \bowtie I$ a proper submodule of $M \bowtie I$. Then $N \bowtie I$ is a prime submodule if and only if $N \bowtie I$ is primary and weakly prime.
\end{corollary}

\begin{proof}
$(\Rightarrow):$
It is clear.

$(\Leftarrow):$ Assume that $N \bowtie I$ is primary and weakly prime submodule, and $(a,a+i)(x,x') \in N \bowtie I$ where $(x,x') \notin N \bowtie I$. Since $N \bowtie I$ is a primary, there is a $n \in \mathbb{Z}^{+}$ such that $(a,a+i)^{n}M \bowtie I \subseteq N \bowtie I$. Therefore, for every $(y,y') \in M \bowtie I \setminus N \bowtie I$, there is a $n \in \mathbb{Z}^{+}$ such that $(a,a+i)^n(y,y') \in N \bowtie I$, i.e., $(a,a+i)^n \in (N \bowtie I:(y,y'))$. Since $N \bowtie I$ is a weakly prime, according to the Lemma \ref{Lemma 3.}, $(i)$, $(N \bowtie I:(y,y'))$ is a prime ideal of $A \bowtie I$. Thus, $(a,a+i) \in (N \bowtie I:(y,y'))$. Hence for all $(y,y') \in M \bowtie I$, we have, $(a,a+i)(y,y') \in N \bowtie I$, that is, $(a,a+i)M \bowtie I \subseteq N \bowtie I$. Therefore, $N \bowtie I$ is a prime submodule.

\end{proof}

\begin{theorem}
\label{Teorem 4.}
    Let $M \bowtie I$ be an $A \bowtie I$-module and $N \bowtie I$ be a proper submodule of $M \bowtie I$. The followings
are equivalent.
\begin{enumerate}
    \item $N\bowtie I$ is a weakly prime submodule.
    \item For any $(x,x'), (y,y') \in M \bowtie I$, if $(N \bowtie I : (x,x')) \neq (N \bowtie I : (y,y'))$, then $N \bowtie I = (N \bowtie I + A \bowtie I(x,x')) \cap (N \bowtie I + A \bowtie I(y,y'))$.
\end{enumerate}

\end{theorem}

\begin{proof}
    $(1) \Rightarrow (2)$ 
    Suppose that $(a,a+i) \in (N \bowtie I: (x,x')) \setminus (N \bowtie I : (y,y'))$, where $(a,a+i) \in A \bowtie I$, i.e., we have $(a,a+i)(x,x') \in N \bowtie I$ and $(a,a+i)(y,y') \notin N \bowtie I$. Hence $(y,y') \notin N \bowtie I$, since $N \bowtie I$ is a weakly prime submodule, according to Lemma \ref{Lemma 3.} $(i)$, $(N \bowtie I : (y,y'))$ is a prime ideal of $A \bowtie I$.\\
    
    We have $(N \bowtie I: (y,y')) \subseteq (N \bowtie I : (a,a+i)(y,y'))$. On the other hand, for all $(b,b+i_1) \in (N \bowtie I : (a,a+i)(y,y'))$, \\ $(b,b+i_1)(a,a+i)(y,y') \in N \bowtie I$ and $(a,a+i)(y,y') \notin N \bowtie I$, that is $(b,b+i_1)(a,a+i) \in (N \bowtie I: (y,y'))$ and $(a,a+i) \notin (N \bowtie I: (y,y'))$. \\
    Since $(N \bowtie I: (y,y'))$ is a prime ideal, then $(b,b+i_1) \in (N \bowtie I : (y,y'))$.\\
    $(N \bowtie I : (a,a+i)(y,y')) \subseteq (N \bowtie I : (y,y'))$, that is, $(N \bowtie I: (y,y')) = (N \bowtie I : (a,a+i)(y,y'))$ is obtained.\\
    If $(t,t+k) \in (N \bowtie I + A \bowtie I(x,x')) \cap (N \bowtie I + A \bowtie I(y,y'))$, then there exists $(a_1,a_1+i_1),(a_2,a_2+i_2) \in A \bowtie I$ and $(n_1, n_1^{'}), (n_2, n_2^{'}) \in N \bowtie I$ such that $(t,t+k) = (n_1, n_1^{'}) + (a_1,a_1+i_1)(x,x') = (n_2, n_2^{'}) + (a_2,a_2+i_2)(y,y')$.\\
   
    We conclude that $(a,a+i)(t,t+k) = (a,a+i)(n_1, n_1^{'}) + (a_1,a_1+i_1)(a,a+i)(x,x') = (a,a+i)(n_2, n_2^{'}) + (a_2,a_2+i_2)(a,a+i)(y,y')$. \\
    Since $(a,a+i)(n_1, n_1^{'}), (a_1,a_1+i_1)(a,a+i)(x,x'),(a,a+i)(n_2, n_2^{'}) \in N \bowtie I $, then we have $(a_2,a_2+i_2)(a,a+i)(y,y') \in N \bowtie I$. That is, $(a_2,a_2+i_2) \in (N \bowtie I: (a,a+i)(y,y')) = ( N \bowtie I: (y,y'))$ is obtained and here, $(a_2,a_2+i_2)(y,y') \in N \bowtie I$. Therefore, we have $(t,t+k) =(n_2, n_2^{'}) + (a_2,a_2+i_2)(y,y') \in N \bowtie I $.\\

    $(2) \Rightarrow (1)$  For $(a_1,a_1+i_1), (a_2,a_2+i_2) \in A \bowtie I$, $(m,m') \in M \bowtie I$, if $(a_1,a_1+i_1)(a_2,a_2+i_2)(m,m') \in N \bowtie I$ and $(a_1,a_1+i_1)(m,m') \notin N \bowtie I$, then it is enough to show that $(a_2,a_2 +i_2)(m,m') \in N \bowtie I$. Hence we can get $(a_1,a_1+i_1) \in (N \bowtie I : (a_2,a_2+i_2)(m,m')) \setminus (N \bowtie I: (m,m'))$. Therefore, we conclude $(N \bowtie I : (a_2,a_2+i_2)(m,m')) \neq (N \bowtie I: (m,m'))$. Now if we take $(x,x') = (a_2,a_2+i_2)(m,m')$, $(y,y') = (m,m')$, from our assumption we have
    $N \bowtie I = (N \bowtie I + A \bowtie I(a_2,a_2+i_2)(m,m')) \cap (N \bowtie I + A \bowtie I(m,m'))$. This shows that $(a_2,a_2+i_2)(m,m') \in (N \bowtie I + A \bowtie I(a_2,a_2+i_2)(m,m')) \cap (N \bowtie I + A \bowtie I(m,m')) = N \bowtie I$.
    
\end{proof}

\begin{remark}
    Let $N \bowtie I$ be a prime submodule in $M \bowtie I$. For $(x,x'), (y,y') \in  M \bowtie I$ and $(a,a+i) \in A \bowtie I$, if $(a,a+i)(x,x') \in N \bowtie I$, then $N \bowtie I = (N \bowtie I)+ A \bowtie I(x,x')$ or $ N \bowtie I = (N \bowtie I) + A \bowtie I(a,a+i)(y,y')$. If $(x,x') \in N \bowtie I$, then we have $A \bowtie I(x,x') \subseteq N \bowtie I$. If $(x,x') \notin N \bowtie I$, then we conclude $(a,a+i)M \bowtie I \subseteq N \bowtie I$, that is $(a,a+i)(y,y') \in N \bowtie I$ and $A \bowtie I(a,a+i)(y,y') \subseteq N \bowtie I$.
\end{remark}

\begin{corollary}
    Let $M \bowtie I$ be an $A \bowtie I$-module, $N \bowtie I$ be a weakly prime submodule of $M \bowtie I$ and $(x,x'), (y,y') \in M \bowtie I$. Then the followings hold.
\begin{enumerate}
    \item  If $(a,a+i)(x,x') \in N \bowtie I$ where $(a,a+i) \in A \bowtie I$, then $N \bowtie I = (N \bowtie I + A \bowtie I(x,x')) \cap (N \bowtie I + A \bowtie I(a,a+i)(y,y'))$.
    \item If $N \bowtie I$ is an irreducible submodule, then $N \bowtie I$ is a prime submodule.
\end{enumerate}
\end{corollary}

\begin{proof}
    \begin{enumerate}
        \item If $(a,a+i)(y,y') \in N \bowtie I$, then $N \bowtie I = (N \bowtie I + A \bowtie I(x,x')) \cap (N \bowtie I + A \bowtie I(a,a+i)(y,y'))$. Assume that $(a,a+i)(y,y') \notin N \bowtie I$. Then, since $(a,a+i) \in (N \bowtie I : (x,x'))$ and $(a,a+i) \notin (N \bowtie I : (y,y'))$, $(N \bowtie I : (x,x')) \neq (N \bowtie I : (y,y'))$. According to Theorem \ref{Teorem 4.}, we have \\ $N \bowtie I \subseteq (N \bowtie I + A \bowtie I(x,x')) \cap (N \bowtie I + A \bowtie I(a,a+i)(y,y')) \subseteq (N \bowtie I + A \bowtie I(x,x')) \cap (N \bowtie I + A \bowtie I(y,y')) = N \bowtie I$.
        \item Suppose that $(a,a+i)(x,x') \in N \bowtie I$, where $(a,a+i) \in A \bowtie I$. According to (1), for every $(y,y') \in M \bowtie I$,\\
$(N \bowtie I + A \bowtie I(x,x')) \cap (N \bowtie I + A \bowtie I(a,a+i)(y,y')) = N \bowtie I$. Since $N \bowtie I$ is an irreducible submodule, then we conclude either \\ $N \bowtie I = (N \bowtie I + A \bowtie I(x,x')) $ or $N \bowtie I = (N \bowtie I + A \bowtie I(a,a+i)(y,y'))$. Therefore $(x,x') \in N \bowtie I$ or $(a,a+i)(y,y') \in N \bowtie I$ is obtained. Hence $N \bowtie I$ is a prime submodule.
    \end{enumerate}
\end{proof}

Let $M$ be an $A$-module and $N$ a submodule of $M$. In this section for every $a \in A$, we consider $(N : a)$ to be:
$(N : a) = \{m \in M| am \in N\}$.
It is easy to see that $(N : a)$ is a submodule of $M$ containing $N$.
We will examine this definition in relation to the duplication of the $A$-module $M$ along the ideal definition.

The following lemma will give us a characterization of weakly prime submodules.

\begin{lemma}
    Let $M \bowtie I$ be an $A \bowtie I$-module and $N \bowtie I$ a proper submodule of $M \bowtie I$. Then $N \bowtie I$ is a weakly prime submodule of $M \bowtie I$ if and only if for every $(a_1,a_1 + i_1), (a_2, a_2 +i_2) \in A \bowtie I$,\\
    $(N \bowtie I : (a_1,a_1 + i_1)(a_2, a_2 +i_2)) = (N \bowtie I: (a_1,a_1 + i_1))$ or \\
    $(N \bowtie I : (a_1,a_1 + i_1)(a_2, a_2 +i_2)) = (N \bowtie I : (a_2, a_2 +i_2))$.
\end{lemma}

\begin{proof}
    $(\Rightarrow):$ Assume that $N \bowtie I$ is a weakly prime submodule of $M \bowtie I$. Firstly, let us show that \\  $(N \bowtie I : (a_1,a_1 + i_1)(a_2, a_2 +i_2)) = (N \bowtie I: (a_1,a_1 + i_1)) \cup (N \bowtie I : (a_2, a_2 +i_2)) $. For every $(x,x') \in (N \bowtie I: (a_1,a_1 + i_1)(a_2, a_2 +i_2))$, \\
    $(a_1,a_1 + i_1)(a_2, a_2 +i_2)(x,x') \in N \bowtie I$, i.e.,  $(a_1,a_1 + i_1)(a_2, a_2 +i_2)(N \bowtie I :(a_1,a_1 + i_1)(a_2, a_2 +i_2)) \subseteq N \bowtie I$. \\
    Since $N \bowtie I$ is a weakly prime submodule, we have either $(a_1,a_1 + i_1)(N \bowtie I :(a_1,a_1 + i_1)(a_2, a_2 +i_2)) \subseteq N \bowtie I$ or $(a_2, a_2 +i_2)(N \bowtie I :(a_1,a_1 + i_1)(a_2, a_2 +i_2)) \subseteq N \bowtie I$. \\
    Then we obtain $(a_1,a_1 + i_1)(x,x') \in N \bowtie I$ or $(a_2, a_2 +i_2)(x,x')  \in N \bowtie I$. In that case, $(x,x') \in (N \bowtie I :(a_1,a_1 + i_1)) \cup (N \bowtie I :(a_2,a_2 + i_2))$.\\
    
    On the other hand, for all $(x,x') \in (N \bowtie I :(a_1,a_1 + i_1)) \cup (N \bowtie I :(a_2,a_2 + i_2))$, we conclude either $(a_1,a_1 + i_1)(x,x') \in N \bowtie I$ or $(a_2, a_2 +i_2)(x,x')  \in N \bowtie I$. Therefore, we have $(a_1,a_1 + i_1)(a_2, a_2 +i_2)(x,x') \in N \bowtie I$. Hence, $(x,x') \in (N \bowtie I: (a_1,a_1 + i_1)(a_2, a_2 +i_2))$ and equality is shown.\\
    
    Assume that $(N \bowtie I: (a_1,a_1 + i_1)(a_2, a_2 +i_2)) \nsubseteq (N \bowtie I: (a_1,a_1 + i_1))$ and $(N \bowtie I: (a_1,a_1 + i_1)(a_2, a_2 +i_2)) \nsubseteq (N \bowtie I: (a_2,a_2 + i_2))$. Then the elements exist $(x,x') \in (N \bowtie I: (a_1,a_1 + i_1)(a_2, a_2 +i_2)) \setminus (N \bowtie I: (a_1,a_1 + i_1)$ and $(y,y') \in (N \bowtie I: (a_1,a_1 + i_1)(a_2, a_2 +i_2)) \setminus (N \bowtie I:(a_2, a_2 +i_2))$. \\ 
    
    Then, according to the equality, $(x,x') + (y,y') \in (N \bowtie I :  (a_1,a_1 + i_1)(a_2, a_2 +i_2)) \subseteq (N \bowtie I : (a_1,a_1 + i_1)) \cup (N \bowtie I : (a_2, a_2 +i_2))$. \\  
    Here,  $(x,x') + (y,y') \in (N \bowtie I : (a_1,a_1 + i_1)) $ or  $(x,x') + (y,y') \in  (N \bowtie I : (a_2, a_2 +i_2))$.\\
    If  $(x,x') + (y,y') \in (N \bowtie I : (a_1,a_1 + i_1)) $, since $(y,y') \in (N \bowtie I: (a_1,a_1 + i_1)(a_2, a_2 +i_2)) \setminus (N \bowtie I:(a_2, a_2 +i_2))$, we have $(y,y') \in (N \bowtie I: (a_1,a_1 + i_1))$ and we have the contradiction $ (x,x') = ((x,x') + (y,y')) - (y,y') \in (N \bowtie I : (a_1,a_1 + i_1))$.\\
    Similarly, if $(x,x') + (y,y') \in (N \bowtie I : (a_2, a_2 +i_2))$, we have the contradiction $(y,y') \in (N \bowtie I: (a_2, a_2 +i_2))$. That is,  $(N \bowtie I: (a_1,a_1 + i_1)(a_2, a_2 +i_2)) \subseteq (N \bowtie I: (a_1,a_1 + i_1))$ or $(N \bowtie I: (a_1,a_1 + i_1)(a_2, a_2 +i_2)) \subseteq (N \bowtie I: (a_2,a_2 + i_2))$.\\

    Furthermore, since $ (N \bowtie I: (a_1, a_1 + i_1)) \subseteq (N \bowtie I : (a_1,a_1 + i_1)) \cup (N \bowtie I : (a_2, a_2 +i_2)) = (N \bowtie I: (a_1,a_1 + i_1)(a_2, a_2 +i_2)) $ and $ (N \bowtie I: (a_2, a_2 + i_2)) \subseteq (N \bowtie I : (a_1,a_1 + i_1)) \cup (N \bowtie I : (a_2, a_2 +i_2)) = (N \bowtie I: (a_1,a_1 + i_1)(a_2, a_2 +i_2)) $, we have  $(N \bowtie I: (a_1,a_1 + i_1)(a_2, a_2 +i_2)) = (N \bowtie I: (a_1,a_1 + i_1))$ or $(N \bowtie I: (a_1,a_1 + i_1)(a_2, a_2 +i_2)) = (N \bowtie I: (a_2,a_2 + i_2))$.\\

    $(\Leftarrow):$ Let $(a_1,a_1 + i_1)(a_2, a_2 +i_2)(x,x') \in N \bowtie I$ where $(a_1,a_1 + i_1), (a_2, a_2 +i_2) \in A \bowtie I$ and $(x,x') \in M \bowtie I$. Suppose that  $(N \bowtie I : (a_1,a_1 + i_1)(a_2, a_2 +i_2)) = (N \bowtie I: (a_1,a_1 + i_1))$. Since $(x,x') \in (N \bowtie I : (a_1,a_1 + i_1)(a_2, a_2 +i_2)) = (N \bowtie I: (a_1,a_1 + i_1)) $, then $(a_1,a_1 + i_1)(x,x') \in N \bowtie I$. So $N \bowtie I$ is a weakly prime submodule. \\
    Similarly, if we take  $(N \bowtie I : (a_1,a_1 + i_1)(a_2, a_2 +i_2)) = (N \bowtie I: (a_2,a_2 + i_2))$, then $(x,x') \in (N \bowtie I : (a_1,a_1 + i_1)(a_2, a_2 +i_2)) = (N \bowtie I: (a_2,a_2 + i_2)) $, then $(a_2,a_2 + i_2)(x,x') \in N \bowtie I$ and $N \bowtie I$ is a weakly prime submodule.
\end{proof}

\begin{remark}
    If $M \bowtie I$ is an $A \bowtie I$-module and $N \bowtie I$ is a prime submodule, then $(N \bowtie I : M \bowtie I)$ is a prime ideal of $A \bowtie I$. However, the same situation is not true for weakly prime submodules. That is, if $N \bowtie I$ is a weakly prime submodule that is not a prime, $(N \bowtie I : M \bowtie I)$ need not be a weakly prime ideal of $A \bowtie I$.  Let consider $A = \mathbb{Z}_6$, $\mathbb{Z}_6$-module $M = \mathbb{Z}_6$ and $N = \Bar{0}$ be submodule of $M$. For $I = 3\mathbb{Z}_6$, $N \bowtie I = \Bar{0} \bowtie 3\mathbb{Z}_6$ is a weakly prime submodule but $(N \bowtie I : M \bowtie I ) = ( \Bar{0} \bowtie 3\mathbb{Z}_6 : \mathbb{Z}_6 \bowtie 3\mathbb{Z}_6)$ is not a weakly prime ideal. For instance, $(\Bar{0}, \Bar{0}) \neq (\Bar{2}, \Bar{5})(\Bar{3}, \Bar{3}) = (\Bar{0}, \Bar{3}) \in (N \bowtie I : M \bowtie I )$ but $(\Bar{2}, \Bar{5}) \notin (N \bowtie I : M \bowtie I )$ and $ (\Bar{3}, \Bar{3}) \notin (N \bowtie I : M \bowtie I )$. However, the following proposition can be given. 
   
\end{remark}

\begin{proposition}
     Let $A \bowtie I$ be a commutative ring, $M \bowtie I$ a faithful cyclic $A \bowtie I$-module, and $N \bowtie I$ a weakly prime submodule of $M 
 \bowtie I$. Then $(N \bowtie I : M \bowtie I)$ is a weakly prime ideal of $A \bowtie I$.
\end{proposition}

\begin{proof}
    Suppose that $M \bowtie I = A \bowtie I(m,m')$ and let $(0,0) \neq (a,a+i)(b,b+j) \in (N \bowtie I : M \bowtie I) $ with $(a,a+i) \notin (N \bowtie I : M \bowtie I)$. Then there exists $(c,c+k) \in A \bowtie I$ such that $(a,a+i)((c,c+k)(m,m')) \notin N \bowtie I$, so we have $(a,a+i)(m,m') \notin N \bowtie I$. Since $(0,0) \neq (a,a+i)(b,b+j)M \bowtie I \subseteq N \bowtie I$, then $(a,a+i)(b,b+j)(m,m') \in N \bowtie I$. \\
    If $(a,a+i)(b,b+j)(m,m') = (0,0)$, since $M \bowtie I$ is a faithful cyclic $A \bowtie I$-module, then we have $(a,a+i)(b,b+j) \in ((0,0): (m,m')) = (0 : M \bowtie I) = (0,0)$ and the contradiction $(a,a+i)(b,b+j) = (0,0)$ would be obtained. Therefore, since $N \bowtie I$ is a weakly prime submodule of $M \bowtie I$, $(b,b+j) \in (N \bowtie I: M \bowtie I)$ is obtained, and so $(N \bowtie I: M \bowtie I)$ is a weakly prime ideal of $A \bowtie I$. \\
\end{proof}

\begin{lemma}
    Let $N \bowtie I$ be a proper submodule of an $A \bowtie I$-module $M \bowtie I$. Then the followings are equivalent:
    \begin{enumerate}
        \item $N \bowtie I$ is a primary submodule of $M \bowtie I$.
        \item For every $ (b,b') \in M \bowtie I \setminus N \bowtie I$ , $ (N \bowtie I : A \bowtie I(b,b')) \subseteq \sqrt{(N \bowtie I: M \bowtie I)}$
    \end{enumerate} 

\end{lemma}

\begin{proof}
   $(1) \Rightarrow (2)$ Suppose that $N \bowtie I$ is a primary submoule of $M \bowtie I$. Let $(c,c+k) \in (N \bowtie I : A \bowtie I(b,b'))$, then we have $(c,c+k)A \bowtie I(b,b') \subseteq N \bowtie I$ for every $(b,b+j) \in M \bowtie I \setminus N \bowtie I$. We conclude that for every $(a,a+i) \in A \bowtie I$, $(c,c+k)(a,a+i)(b,b') \in N \bowtie I$, and so $(c,c+k)(b,b') \in N \bowtie I$.  \\
   $(b,b') \in M \bowtie I - N \bowtie I$ and since $N \bowtie I$ is a primary submodule, then we have $(c,c+k)^n \in (N \bowtie I : M \bowtie I)$. Therefore, $(c,c+k) \in \sqrt{(N \bowtie I : M \bowtie I)}$.\\

   $(2) \Rightarrow (1)$ Suppose that for every $ (b,b') \in M \bowtie I \setminus N \bowtie I$ , $ (N \bowtie I : A \bowtie I(b,b')) \subseteq \sqrt{(N \bowtie I: M \bowtie I)}$. Let $(a,a+i)(b,b') \in N \bowtie I$. Hence for all $(c,c+i_1) \in A \bowtie I$, we get $(a,a+i)(c,c+i_1)(b,b') \in N \bowtie I$. Therefore we have $(a,a+i) \in (N \bowtie I : A \bowtie I(b,b'))$. Since  $ (N \bowtie I : A \bowtie I(b,b')) \subseteq \sqrt{(N \bowtie I: M \bowtie I)}$, then $(a,a+i) \in \sqrt{(N \bowtie I: M \bowtie I)}$. So, there exist an $n \in \mathbb{Z}^{+}$ such that $(a,a+i)^n \in (N \bowtie I: M \bowtie I)$.

\end{proof}

\begin{proposition}
   Let $M \bowtie I$ be a $A \bowtie I$-module and let $N \bowtie I$ be a primary submodule in $M \bowtie I$. Then $(N \bowtie I: M \bowtie I) = Ann(M \bowtie I/N \bowtie I)$ is the primary ideal in $A \bowtie I$.
\end{proposition}

\begin{proof}
    Let $N \bowtie I$ is the primary submodule of $M \bowtie I$. Assume that $(a_1,a_1+i_1)(a_2,a_2+i_2) \in (N \bowtie I : M \bowtie I)$ and $(a_2,a_2+i_2) \notin (N \bowtie I : M \bowtie I)$ where $(a_1,a_1+i_1), (a_2,a_2+i_2) \in A \bowtie I$. Hence, $(a_1,a_1+i_1)(a_2,a_2+i_2)M \bowtie I \subseteq N \bowtie I$ and $(a_2,a_2 +i_2)M \bowtie I \nsubseteq N \bowtie I$, that is, there is at least one $(m,m') \in  M \bowtie I$ such that $(a_1,a_1+i_1)(a_2,a_2 + i_2)(m,m') \in N \bowtie I$ and $(a_2,a_2+i_2)(m,m') \notin N \bowtie I$. Since $N \bowtie I$ is primary, there exists a positive integer $c$ such that $(a_1,a_1 +i_1)^cM \bowtie I \subseteq N \bowtie I$. Therefore, $(a_1,a_1 +i_1)^c \in (N \bowtie I: M \bowtie I)$. Thus,  $(N \bowtie I: M \bowtie I)$ is a primary ideal of $A \bowtie I$.
\end{proof}

\begin{corollary}
    Let $M \bowtie I$ be $A \bowtie I$-module, $N \bowtie I$ be the primary submodule of $M \bowtie I$. $\sqrt{(N \bowtie I: M \bowtie I)}$ is a prime ideal in $A \bowtie I$.
\end{corollary}

\begin{lemma}
\label{Lemma 8.}
    $0 \times IM$ and $IM \times IM$ are $A \bowtie I$-submodules of $M \bowtie I$. The following isomorphisms hold \cite{bouba2018duplication}:
   \begin{center}
       \begin{enumerate}
           \item $\frac{M \bowtie I}{0 \times IM} \cong M$
           \item $\frac{M \bowtie I}{IM \times IM} \cong \frac{M}{IM}$
       \end{enumerate}
   \end{center} 
\end{lemma}

\noindent A left $A$-module $M$ is called weakly prime if the annihilator of any nonzero submodule of $M$ is a prime ideal and a proper submodule $P \subset M$ is called weakly prime submodule if the quotient module $M / P$ is a weakly prime module \cite{behboodi2004weakly}.
\,\\

\begin{theorem} 
\item[(1)]  $M \bowtie I$ is a weakly prime $A \bowtie I$-module if and only if $IM = 0$ and $M$ is a weakly prime $A$-module.
\item[(2)] $0 \times IM$ is a weakly prime $A \bowtie I$-submodule of $M \bowtie I$ if and only if $M$ is a weakly prime $A$-module.
\end{theorem}

\begin{proof}
\begin{enumerate}
    \item $(\Rightarrow):$ Let $M \bowtie I$ be weakly prime $A \bowtie I$-module and suppose that $IM \neq 0$. Thus, there exist $i \in I$ and $m \in M$ such that $im \neq 0$. From here we have $(i, 0)(0, im) = (0, 0)$ while $(0, im) \neq (0,0)$. Since $(0,i)(i,m)=(0,im) \in Ann_{l}(M \bowtie I)$ and $M \bowtie I$ is a weakly prime module, it must be $(0,i) \in Ann_{l}(M \bowtie I)$ or $(i,m) \in Ann_{l}(M \bowtie I)$. \\
    If $(0,i) \in Ann_{l}(M \bowtie I)$, for all $(m',m) \in M \bowtie I$, it must be $(0,i)(m',m) = (0,im) = (0,0)$ which contradicts the fact that $im \neq 0$. Then, $IM = 0$. \\
    Also if $(i,m) \in Ann_{l}(M \bowtie I)$, for all $(m, m') \in M \bowtie I$, it must be $(i,m)(m,m') = (im,mm') = (0,0)$ which contradicts the fact that $im \neq 0$. Thus we obtain $IM = 0$.
    Now, let $a_1a_2m = 0$ where $ 0 \neq a_1, a_2 \in A, m \in M$. Therefore, we have $a_1a_2 \in Ann_{l}(M)$. 
    Now, $(a_1a_2, a_1a_2)(m, m) = (0,0)$ and thus $(a_1,a_1)(a_2,a_2) = (a_1a_2, a_1a_2) \in Ann_{l}(M \bowtie I) $. Since $M \bowtie I$ is a weakly prime module, then it must be $(a_1,a_1) \in Ann_{l}(M \bowtie I)$ or $(a_2, a_2) \in Ann_{l}(M \bowtie I)$. \\
    If $(a_1,a_1) \in Ann_{l}(M \bowtie I)$, then for all $(m_1,m_1) \in M \bowtie I$, we have $(a_1,a_1)(m_1,m_1) = (0,0)$. Hence, we obtain $a_1m_1 = 0$.\\
    If $(a_2, a_2) \in Ann_{l}(M \bowtie I)$, then it must be $(a_2,a_2)(m_2,m_2) = (0,0)$ for all $(m_2,m_2) \in M \bowtie I$, and we obtain $a_2m_2 = 0$. \\
    We conclude that $a_1 \in Ann_{l}(M)$ or $a_2 \in Ann_{l}(M)$ and thus $Ann_{l}(M)$ is a prime ideal. Therefore, $M$ is a weakly prime $A$-module.\\
    
    $(\Leftarrow):$ Let $IM = 0$ and $M$ be a weakly prime $A$-module. Since $IM = 0$, we have $M \bowtie I = \{(m, m)|m \in M\}$. Now, assume that $(a_1,a_1+i_1)(a_2, a_2 + i_2)(m,m) = (0,0)$ where $(0,0) \neq (a_1,a_1+i_1),(a_2, a_2 + i_2) \in A \bowtie I$ and $(m,m) \in M \bowtie I$. Thus, we have $(a_1,a_1+i_1)(a_2, a_2 + i_2) \in Ann_{l}(M \bowtie I)$. Therefore, $a_1a_2m = 0$, $a_1a_2 \in Ann(M)$ and $i_1,i_2 \in Ann_{l}(M)$. Since $M$ is a weakly prime module, then it must be $a_1 \in Ann_{l}(M)$ or $a_2 \in Ann_{l}(M)$. \\
    If $a_1 \in Ann_{l}(M)$, then for all $m_1 \in M$, it must be $a_1m_1 = 0$. Hence, we obtain $(a_1,a_1 + i_1)(m_1,m_1) =(0,0)$.  \\
    Also if $a_2 \in Ann_{l}(M)$, then for all $m_2 \in M$, it must be $a_2m_2 = 0$ and we have $(a_2,a_2 +i_2)(m_2,m_2) =(0,0)$. \\
    Consequently, we conclude that $(a_1,a_1 + i_1) \in Ann_{l}(M \bowtie I)$ or $(a_2,a_2 + i_2) \in  Ann_{l}(M \bowtie I)$. Therefore, $Ann_{l}(M \bowtie I)$ is a prime ideal and so $M \bowtie I$ is a weakly prime $A \bowtie I$-module. \\
    
    \item 
    According to the Lemma \ref{Lemma 8.}, if $M$ is a weakly prime $A$-module, then $\frac{M \bowtie I}{0 \times IM} \cong M$ is a weakly prime module. By definition, $0 \times IM$ is a weakly prime $A \bowtie I$-submodule of $M \bowtie I$. Conversely, the desired result follows from Lemma \ref{Lemma 8.}.
 \end{enumerate}
    
\end{proof}

$\,$

$\,$


\begin{thebibliography}{99}

\bibitem{atani2007weakly} S. Ebrahimi Atani and F. Farzalipour. $On$ $weakly$ $prime$ $submodules$.
Tamkang Journal of Mathematics, 38(3):247–252, 2007. 
	\bibitem{d2007amalgamated} Marco D’Anna and Marco Fontana. $\textit{An amalgamated duplication of a ring along an ideal:} \\ \textit{the basic properties}$. Journal of Algebra and its Applications, 6(03):443–459, 2007.
	\bibitem{issoual2022s} Mohammed Issoual, Najib Mahdou, Neslihan Aysen Ozkirisci, and
Ece Yetkin Celikel. \textit{S-1-absorbing primary submodules.} arXiv preprint
arXiv:2203.04690, 2022.
	\bibitem{bouba2018duplication} EM Bouba, N Mahdou, and M Tamekkante. \textit{ Duplication of a module along
an ideal}. Acta Mathematica Hungarica, 154:29–42, 2018.
	\bibitem{behboodi2004weakly} M Behboodi and H Koohy. \textit{Weakly prime modules}. Vietnam J. Math,
32(2):185–195, 2004.

\end{thebibliography}
\end{document}